\documentclass{article}

\usepackage{amsthm}
\usepackage{amsmath}
\usepackage{amssymb}
\usepackage{dsfont}
\usepackage[all]{xy}
\usepackage{tikz-cd}
\usepackage{enumitem}
\usepackage{newclude}
\usepackage{hyperref}

\theoremstyle{definition}
\newtheorem{Definition}{Definition}[section]

\theoremstyle{plain}
\newtheorem{Theorem}[Definition]{Theorem}

\theoremstyle{plain}
\newtheorem{Proposition}[Definition]{Proposition}

\theoremstyle{plain}
\newtheorem{Lemma}[Definition]{Lemma}

\theoremstyle{plain}
\newtheorem{Corollary}[Definition]{Corollary}

\theoremstyle{definition}
\newtheorem{Example}[Definition]{Example}

\theoremstyle{definition}
\newtheorem{Notation}[Definition]{Notation}

\theoremstyle{definition}

\theoremstyle{definition}

\theoremstyle{remark}
\newtheorem{Remark}[Definition]{Remark}

\title{The Relative Deligne Tensor Product over Pointed Braided Fusion Categories}
\author{Thibault D. Décoppet}
\date{June 2022}

\AtEndDocument{
    \noindent\textsc{Mathematical Institute, University of Oxford,  Oxford OX2 6GG, United Kingdom}\\
    \textit{E-mail address}: \texttt{thibault.decoppet@maths.ox.ac.uk}\\
    \textit{URL}: \texttt{https://www.thibaultdecoppet.com}\\
}

\begin{document}

\bibliographystyle{alpha}

\maketitle

    \begin{abstract}
        We give a formula for the relative Deligne tensor product of two indecomposable finite semisimple module categories over a pointed braided fusion category over an algebraically closed field.
    \end{abstract}
    
\section*{Introduction}

The relative Deligne tensor product of module categories is a central construction of the theory of finite tensor categories, which categorifies the relative tensor product of modules over an algebra. For instance, working over an algebraically closed field of characteristic zero, it is shown in \cite{ENO2} that, for a fixed fusion category $\mathcal{C}$, graded extensions of $\mathcal{C}$ are classified by maps into the Brauer-Picard 3-groupoid $BrPic(\mathcal{C})$ of invertible $\mathcal{C}$-$\mathcal{C}$-bimodule categories. In this case, the relative Deligne tensor product over $\mathcal{C}$ provides the operation of composition of 1-morphisms. Another example is given by considering $\mathcal{B}$ a braided fusion category. Then, the relative Deligne tensor product over $\mathcal{B}$ of two finite semisimple right $\mathcal{B}$-module categories possesses a right action by $\mathcal{B}$. Thus, we can consider the Picard 3-groupoid $Pic(\mathcal{B})$ of invertible right $\mathcal{B}$-module categories, which was shown to classify crossed braided extensions of $\mathcal{B}$ in \cite{ENO2}. The two preceding examples have been recently shown to hold more generally for any (braided) finite tensor categories over an arbitrary algebraically closed field (see \cite{DN2}). Now, by definition, the Brauer-Picard and Picard groupoids only know about invertible (bi)module categories. But the relative Deligne tensor product provides a natural composition operation for all exact bimodule categories. This can be used to construct a symmetric monoidal 3-category whose objects are finite tensor categories, and 1-morphisms are exact bimodule categories (see \cite{DSPS13}). Via the cobordism hypothesis, the study of this symmetric monoidal 3-category initiated in \cite{DSPS13} has highlighted the link between finite tensor categories and low dimensional topology. In a somehwat different direction, given $\mathcal{C}$ a multifusion category over an algebraically closed field of characteristic zero, we may form $\mathbf{Mod}(\mathcal{C})$ the 2-category of all finite semisimple right $\mathcal{C}$-module categories. The properties of these 2-categories were axiomatized in \cite{DR} under the definition of finite semisimple 2-categories. Further, if $\mathcal{B}$ is a braided fusion category, the 2-category $\mathbf{Mod}(\mathcal{B})$ inherits a rigid monoidal structure given by the relative Deligne tensor product over $\mathcal{B}$, and so is a fusion 2-category (see \cite{DR} for the definition). In particular, the maximal subgroupoid of $\mathbf{Mod}(\mathcal{B})$ is $Pic(\mathcal{B})$. Analogously, for a general fusion category $\mathcal{C}$, we recover $BrPic(\mathcal{C})$ as the maximal subgroupoid of $\mathbf{Mod}(\mathcal{Z}(\mathcal{C}))$, the fusion 2-category associated to the Drinfel'd center of $\mathcal{C}$. These last considerations continue to hold over an arbitrary algebraically closed field as discussed in \cite{D5}.

Due to their connections with graded extensions, the Brauer-Picard group and the Picard group have received a lot of attention (see \cite{ENO2}, \cite{DN1}, \cite{BN}, \cite{Mo}, \cite{GJS}, \cite{GS}, \cite{GP}, \cite{EdM}, \cite{MN}, and \cite{NR}). But, relatively few explicit computations of the relative Deligne tensor product of non-invertible exact (bi)module categories are known. For instance, the relative Deligne tensor product over the symmetric fusion category of representations of a finite group was computed in \cite{Gre}. Another example is given as follows. Let $A$ be a finite abelian group, and write $\mathbf{Vect}_A$ for the fusion category of $A$-graded finite dimensional vector spaces. Working over an algebraically closed field of characteristic zero, a formula for the relative Deligne tensor product of all finite semisimple module categories over $\mathbf{Vect}_A$ with the trivial symmetric braiding was given in \cite{ENO2}. In \cite{D2}, we extended their computations to obtain a formula for the relative Deligne tensor product over $\mathbf{Vect}_A^{\beta}$, where $\beta$ is any braiding on $\mathbf{Vect}_A$. In the present article, we generalize this result in two directions. Firstly, we will allow our base filed to be any algebraically closed field. Secondly, we will let the underlying fusion category carry any associator $\omega$, and compute the relative Deligne tensor product of all finite semisimple module categories over $\mathbf{Vect}_A^{(\omega,\beta)}$, where $\beta$ is a compatible braiding. As a consequence, we get an explicit formula for the fusion rule of $\mathbf{Mod}(\mathbf{Vect}_A^{(\omega,\beta)})$. Finally, we give tables for the relative Deligne tensor product over $\mathbf{Vect}_A^{(\omega,\beta)}$, when $A$ is either $\mathbb{Z}/p$, $\mathbb{Z}/p^2$, $\mathbb{Z}/2\oplus\mathbb{Z}/2$, or $\mathbb{Z}/4\oplus\mathbb{Z}/2$ with various choices of associator and compatible braiding. We also compute the Picard groups of $\mathbb{Z}/p\oplus \mathbb{Z}/p$ and $\mathbb{Z}/2\oplus\mathbb{Z}/2\oplus\mathbb{Z}/2$ with specific choices of associators and compatible braidings. Let us conclude by mentioning that, for some specific groups and braidings, surprisingly similar formulas have been independently obtained in \cite{RSS} through high energy physics considerations. The precise relation with our present work remains to be clarified.

\subsubsection*{Acknowledgments}

I would like to thank Christopher Douglas for pushing me to write this article, as well as for suggesting various improvements.

\section{Preliminaries}

Let $\mathds{k}$ be an algebraically closed field, and $\mathcal{C}$ be a fusion category over $\mathds{k}$. A right $\mathcal{C}$-module category is a category $\mathcal{M}$ over $\mathds{k}$ together with a linear right action $\mathcal{M}\times \mathcal{C}\rightarrow \mathcal{M}$ and natural isomorphisms $$\mu_{M,C,D}:M\otimes (C\otimes D) \rightarrow (M\otimes C)\otimes D,$$ for every $M\in\mathcal{M}$ and $C,D\in\mathcal{C}$ that are suitably coherent. We call such a module category indecomposable if it cannot be written as the direct sum of two non-trivial module subcategories. Every finite semisimple right $\mathcal{C}$-module category can be written as a direct sum of indecomposable finite semisimple right $\mathcal{C}$-module categories (see \cite{EGNO}).

Thanks to a theorem of Ostrik (see \cite{O}, and \cite{DSPS14} for its generalization to arbitrary fields), every finite semisimple right $\mathcal{C}$-module category $\mathcal{M}$ is equivalent to the right $\mathcal{C}$-module category $Mod_{\mathcal{C}}(R)$ of left $R$-modules in $\mathcal{C}$ for some algebra $R$ in $\mathcal{C}$. As this will be important for us later, let us be more precise. One can construct a bilinear functor $\underline{Hom}(-,-):\mathcal{M}^{op}\times\mathcal{M}\rightarrow \mathcal{C}$, which provides an enrichment of $\mathcal{M}$ over $\mathcal{C}$, and fits into an adjunction $$Hom_{\mathcal{M}}(M\otimes C, N)\cong Hom_{\mathcal{C}}(C, \underline{Hom}(M,N)),$$ with $C\in\mathcal{C}$ and $M,N\in\mathcal{M}$. Now, for a sufficiently general choice of $P$ in $\mathcal{M}$, the functor \begin{equation}\label{eqn:adjunction}\underline{Hom}(P,-):\mathcal{M}\rightarrow Mod_{\mathcal{C}}(\underline{End}(P))\end{equation} is an equivalence. If $\mathcal{M}$ is indecomposable, then we may chose $P$ to be any non-zero object.

We will also be interested in the following special class of module categories. We say that a finite semisimple right $\mathcal{C}$-module category $\mathcal{M}\simeq Mod_{\mathcal{C}}(R)$ is separable if $R$ is a separable algebra, i.e. the multiplication map $R\otimes R\rightarrow R$ has a section as a map of $R$-$R$-bimodules. By \cite{DSPS13}, this does not depend on the choice of $R$. Moreover, they show that if $char(\mathds{k})=0$, then every finite semisimple $\mathcal{C}$-module category is separable.

Now, assume that $\mathcal{C}$ admits a braiding $\beta$, and let $\mathcal{M}\simeq Mod_{\mathcal{C}}(R)$ and $\mathcal{N}\simeq Mod_{\mathcal{C}}(S)$ be two finite semisimple right $\mathcal{C}$-module categories. Note that $R\otimes S$ inherits an algebra structure via 
\begin{align}(R\otimes S)\otimes (R\otimes S)&\xrightarrow{\omega} R\otimes (S\otimes (R\otimes S))\xrightarrow{Id\otimes\omega^{-1}}R\otimes ((S\otimes R)\otimes S)\nonumber \\ & \xrightarrow{Id\otimes (\beta\otimes Id)} R\otimes ((R\otimes S)\otimes S)\xrightarrow{Id\otimes \omega} R\otimes (R\otimes (S\otimes S))\nonumber \\ &\xrightarrow{\omega^{-1}}(R\otimes R)\otimes (S\otimes S)\rightarrow R\otimes S,\label{eqn:tensoralgebras}\end{align}
where $\omega$ denotes the associator of $\mathcal{C}$. In particular, we may consider the right $\mathcal{C}$-module category $Mod_{\mathcal{C}}(R\otimes S)$, which is called the relative Deligne tensor product of $\mathcal{M}$ and $\mathcal{N}$ over $\mathcal{C}$, and is denoted by $\mathcal{M}\boxtimes_{\mathcal{C}}\mathcal{N}$. The relative Deligne tensor product can be characterized using a universal property (see \cite{ENO2} for details), thence it does not depend on the choice of algebras $R$ and $S$. Furthermore, under our hypotheses, it is always finite semisimple (see \cite{DSPS14}).

Let us now consider $\mathbf{Mod}(\mathcal{C})$ the 2-category of separable right $\mathcal{C}$-module categories. In the language of \cite{DR} (see also \cite{D5}), this is a finite semisimple 2-category, and its simple objects are precisely the indecomposable separable right $\mathcal{C}$-module categories. Now, it was shown in \cite{DSPS13} that, if $\mathcal{M}$ and $\mathcal{N}$ are separable, then so is $\mathcal{M}\boxtimes_{\mathcal{C}}\mathcal{N}$. In particular, the relative Deligne tensor product endows $\mathbf{Mod}(\mathcal{C})$ with a rigid monoidal structure, so that $\mathbf{Mod}(\mathcal{C})$ is an example of a connected fusion 2-category (see \cite{D5} for more details), and determining its fusion rule amounts to computing the relative Deligne tensor product of the indecomposable separable right $\mathcal{C}$-module categories.

\subsection{Module categories over pointed fusion categories}

Let $G$ be a finite group. We use $\mathbf{Vect}_G$ to denote the (skeletal) category of finite dimensional $G$-graded $\mathds{k}$-vector spaces, and we write $\mathds{k}_g$ with $g\in G$ for the one dimensional vector space with grading $g$. There is a canonical product on this category given by $\mathds{k}_g\otimes \mathds{k}_h = \mathds{k}_{gh}$. Now, given a 3-cocycle $\omega:G\times G\times G\rightarrow \mathds{k}^{\times}$, we get a monoidal structure on $\mathbf{Vect}_G$ with associator given by $$\omega(g,h,m):(\mathds{k}_g\otimes \mathds{k}_h)\otimes \mathds{k}_m\rightarrow \mathds{k}_g\otimes (\mathds{k}_h\otimes \mathds{k}_m),$$ where $g,h,m\in G$. We denote this fusion category by $\mathbf{Vect}^{\omega}_G$. Note that $\mathbf{Vect}^{\omega}_G$ is pointed, i.e. every simple object is invertible.

\begin{Remark}
By theorem 2.11.5 of \cite{EGNO}, every pointed fusion category is equivalent to $\mathbf{Vect}^{\omega}_G$ for some finite group $G$ and 3-cocycle $\omega:G\times G\times G\rightarrow \mathds{k}^{\times}$. Further, cohomologous 3-cocycle yield equivalent fusion categories. In particular, without loss of generality, we may assume that $\omega$ is normalized.
\end{Remark}

\begin{Notation}
Let $H$ be a subgroup of $G$, and $\phi:H\times H\rightarrow\mathds{k}^{\times}$ a 2-cochain such that $d\phi = \omega|_{H}^{-1}$. We write $R(H,\phi)$ for the algebra in $\mathbf{Vect}^{\omega}_G$ given by the twisted group algebra $\mathds{k}^{\phi}\lbrack H\rbrack$ with grading induced by $H\subseteq G$. We denote by $\mathcal{M}(H,\phi)$ the right $\mathbf{Vect}^{\omega}_G$-module category of left $R(H,\phi)$-modules in $\mathbf{Vect}^{\omega}_G$.
\end{Notation}

The following result is standard (see for instance example 7.4.10 in \cite{EGNO}). For later use, we recall the proof given in example 2.6 of \cite{Nai}.

\begin{Lemma}\label{lem:algclosedmodules}
Every indecomposable finite semisimple right $\mathbf{Vect}^{\omega}_G$-module category is equivalent to $\mathcal{M}(H,\phi)$, for some subgroup $H$ of $G$ and 2-cochain $\phi:H\times H\rightarrow\mathds{k}^{\times}$ satisfying $d\phi = \omega|_{H}^{-1}$. 
\end{Lemma}
\begin{proof}
Without loss of generality, we may assume that $\mathcal{M}$ is skeletal. As $\mathcal{M}$ is indecomposable, its set of simple objects is a transitive $G$-set, i.e. it can be written as $H\backslash G$ for some subgroup $H$ of $G$. We write $\mathds{k}_x$ to denote the simple object of $\mathcal{M}$ corresponding to $x\in H\backslash G$. Using the coherence isomorphism of $\mathcal{M}$, we define a 2-cochain $\mu:G\times G\rightarrow Fun(H\backslash G,\mathds{k}^{\times})$ by $\mu(g,h) (x) = \mu_{\mathds{k}_x,\mathds{k}_g,\mathds{k}_h}$ for every $g,h\in G$ and $x\in H\backslash G$. The coherence axiom for right module categories implies that $d\mu = \omega^{-1}$, where $\omega$ is viewed as 3-cocycle whose values are constant functions in $Fun(H\backslash G,\mathds{k}^{\times})$.

Let us now fix $x\in H\backslash G$. Following the proof of Ostrik's theorem, we now wish to determine the algebra $\underline{Hom}(\mathds{k}_x,\mathds{k}_x)$ in $\mathbf{Vect}^{\omega}_G$. By the adjunction (\ref{eqn:adjunction}), we find that the underlying object is $\mathds{k}\lbrack H\rbrack$ with grading induced by $H\subseteq G$. Moreover, tracing through the ajunction carefully, we find that the multiplication on $\mathds{k}\lbrack H\rbrack$ is the usual one twisted by the 2-cochain $\phi:H\times H\rightarrow \mathds{k}^{\times}$ given by $\phi(h_1,h_2) = \mu(h_1,h_2)(x)$ for $h_1,h_2\in H$. In particular, we have $d\phi = \omega|_H$.
\end{proof}

\begin{Remark}\label{rem:Schapirolemma}
Observe that if $\phi$ and $\phi'$ are two 2-cochains $H\times H\rightarrow \mathds{k}^{\times}$ that differ by a coboundary, then $R(H,\phi)\cong R(H,\phi')$ as algebras, and so $\mathcal{M}(H,\phi)\simeq \mathcal{M}(H,\psi)$ as module categories. In fact, more can be said.
\end{Remark}

\begin{Proposition}\label{prop:algclosedmodules}
Two right $\mathbf{Vect}^{\omega}_G$-module categories $\mathcal{M}(H,\phi)$ and $\mathcal{M}(L,\psi)$ are equivalent if and only if there exists an element $g$ in $G$ such that $H = gLg^{-1}$ and $\psi^{-1}\phi^{g}X_g:L\times L\rightarrow \mathds{k}^{\times}$ is a coboundary, where $$\phi^{g}(l,k):=\phi(^gl,{^gk})\ \ \mathrm{and}\ \ X_g(l,k):=\frac{\omega(^gl,g,k)}{\omega(^gl,{^gk},g)\omega(g,l,k)}.$$
\end{Proposition}
\begin{proof}
The proof follows immediately by considering an arbitrary algebraically closed field in \cite{Nat}. Namely, the characteristic zero assumption that is made there is never used. We wish to point out that the fusion category we denote by $\mathbf{Vect}^{\omega}_G$ corresponds to $\mathcal{C}(G,\omega^{-1})$ in the notation of \cite{Nat}. In particular, $R(E,\phi)$ and $R(F,\psi)$ correspond to $A(E,\phi)$ and $A(F,\psi)$ in the notation of \cite{Nat}. Further, the $\Omega_g$ for $\mathcal{C}(G,\omega^{-1})$ appearing in the statement of theorem 1.1 of \cite{Nat} is equal to $X_g^{-1}$ as defined above.
\end{proof}

\begin{Remark}
If we have $H=L$ as in the statement of proposition \ref{prop:algclosedmodules}, then we may take $g=e$, and we find that $X_e$ is a coboundary. So the above proposition is a refinement of remark \ref{rem:Schapirolemma}.
\end{Remark}

\begin{Corollary}\label{cor:abelianmodulecategoriesequivalences}
Let $G=A$ be abelian. Then, the set of equivalence classes of indecomposable finite semisimple right $\mathbf{Vect}^{\omega}_A$-module categories is in bijection with pairs consisting of a subgroup $H\subseteq G$ on which $\omega$ is trivializable together with a class in $H^2(H;\mathds{k}^{\times})$.
\end{Corollary}

It remains to characterize those finite semisimple right $\mathbf{Vect}^{\omega}_G$-module categories that are separable. As we have recalled above, every finite semisimple module category is separable if $char(\mathds{k})=0$, so it is enough to consider the positive characteristic case.

\begin{Lemma}\label{lem:charsep}
If $char(\mathds{k}) = p>0$. The finite semisimple module category $\mathcal{M}(H,\phi)$ is separable if and only if $H$ has no $p$-torsion.
\end{Lemma}
\begin{proof}
Assume $H$ has $p$-torsion. By Cauchy's theorem, there exists an element $h$ of order $p$ in $H$. Let us consider the injective tensor functor $\mathbf{Vect}_{\mathbb{Z}/p}\rightarrow \mathbf{Vect}^{\omega}_G$ corresponding to the group homomorphism $1\mapsto h$. (Such a monoidal functor exists as $\omega$ restricts to coboundary on $H$.) By duality in the form of theorem 7.17.4 of \cite{EGNO}, the tensor functor $$End_{\mathbf{Vect}^{\omega}_G}(\mathcal{M}(H,\phi))\rightarrow End_{\mathbf{Vect}_{\mathbb{Z}/p}}(\mathcal{M}(H,\phi))$$ is surjective. Now, as right $\mathbf{Vect}_{\mathbb{Z}/p}$-module categories, we have $\mathcal{M}(H,\phi)\simeq \oplus_{g\in H\backslash G}\mathbf{Vect}$. Thus, we find $$End_{\mathbf{Vect}_{\mathbb{Z}/p}}(\mathcal{M}(H,\phi))\simeq \bigoplus_{g_1,g_2\in H\backslash G} End_{\mathbf{Vect}_{\mathbb{Z}/p}}(\mathbf{Vect})$$ as a finite category. But, $End_{\mathbf{Vect}_{\mathbb{Z}/p}}(\mathbf{Vect})\simeq \mathrm{Rep}(\mathds{k}\lbrack \mathbb{Z}/p\rbrack)$ is not semisimple, and so appealing to theorem 2.5.4 of \cite{DSPS13} proves that $\mathcal{M}(H,\phi)$ is not a separable $\mathbf{Vect}^{\omega}_G$-module category.

Now, if we assume that $H$ has no $p$-torsion, then the algebra $R(H,\phi)=\mathds{k}^{\phi}\lbrack H\rbrack$ in $\mathbf{Vect}^{\omega}_G$ is separable. Namely, using the equality $d\phi = \omega|_H^{-1}$, we find that the map $\Delta:\mathds{k}^{\phi}\lbrack H\rbrack\rightarrow \mathds{k}^{\phi}\lbrack H\rbrack\otimes \mathds{k}^{\phi}\lbrack H\rbrack$ given by $$h\mapsto \frac{1}{|H|}\sum_{f\in H}\frac{1}{\phi(f^{-1},fh)}f^{-1}\otimes fh,$$ is a section of the multiplication map of $\mathds{k}^{\phi}\lbrack H\rbrack$, and is a map of $\mathds{k}^{\phi}\lbrack H\rbrack$-$\mathds{k}^{\phi}\lbrack H\rbrack$-bimodules. This concludes the proof of the lemma.
\end{proof}

\subsection{Pointed braided fusion categories}

Let $A$ be a finite abelian group. Recall from \cite{ML} that an abelian 3-cocycle for $A$ with coefficient in $\mathds{k}^{\times}$ is a pair $(\omega, \beta)$ consisting of a 3-cocycle $\omega:A\times A\times A\rightarrow \mathds{k}^{\times}$ and a function $\beta:A\times A\rightarrow \mathds{k}^{\times}$ satisfying $$\omega(b,c,a)\beta(a,b+c)\omega(a,b,c) = \beta(a,c)\omega(b,a,c)\beta(a,b), $$ $$\omega(c,a,b)^{-1}\beta(a+b,c)\omega(a,b,c)^{-1}=\beta(a,c)\omega(a,c,b)^{-1}\beta(b,c),$$ for every $a,b,c\in A$. An abelian 3-cocycle for $A$ with coefficient in $\mathds{k}^{\times}$ is a coboundary if there exists a function $\kappa:A\times A\rightarrow \mathds{k}^{\times}$ such that $$\omega(a,b,c) = \kappa(b,c)\kappa(a+b,c)^{-1}\kappa(a,b+c)\kappa(a,b)^{-1},$$ $$\beta(a,b) = \kappa(a,b)\kappa(b,a)^{-1},$$ for every $a,b,c\in A$. We write $H^3_{ab}(A;\mathds{k}^{\times})$ for the quotient of the abelian group of abelian 3-cocycle by the subgroup of coboundaries. It was proven in \cite{ML} that $H^3_{ab}(A;\mathds{k}^{\times})$ is isomorphic to the group of quadratic forms $q:A\rightarrow \mathds{k}^{\times}$. The isomorphism is given by sending $(\omega,\beta)$ to the quadratic form $a\mapsto \beta(a,a)$.

Now, let $A$ be a finite abelian group and $(\omega,\beta)$ an abelian 3-cocycle for $A$ with coefficient in $\mathds{k}^{\times}$. We write $\mathbf{Vect}^{(\omega,\beta)}_A$ for the pointed fusion category $\mathbf{Vect}^{\omega}_A$ equipped with the braiding given by $$\beta(a,b):\mathds{k}_a\otimes \mathds{k}_b\rightarrow \mathds{k}_b\otimes \mathds{k}_a,$$ for every $a,b\in A$. Given $q:A\rightarrow \mathds{k}^{\times}$ a quadratic form, Quinn's formula (see subsection 2.5.2 of \cite{Qui}) gives an explicit abelian 3-cocycle corresponding to $q$, and we write $\mathbf{Vect}_A^q$ for the associated pointed braided fusion category.

\begin{Remark}\label{rem:pointedbraided}
By proposition 3.1 of \cite{JS}, every pointed braided fusion category is equivalent to $\mathbf{Vect}^{(\omega,\beta)}_A$ for some finite abelian group $A$ and abelian 3-cocycle $(\omega,\beta)$. Further, equivalent abelian 3-cocycles give equivalent pointed braided fusion categories.
\end{Remark}

The following lemma and its corollary are useful when describing the indecomposable finite semisimple module categories over $\mathbf{Vect}^{(\omega,\beta)}_A$.

\begin{Lemma}\label{lem:subgrouptrivializable}
Let $q(a)=\beta(a,a)$ for all $a\in A$ be the quadratic form associated to the abelian 3-cocycle $(\omega,\beta)$. The restriction of the 3-cocycle $\omega$ to the subgroup $E$ of $A$ is trivializable if and only if, for every $e\in E$, the order of $q(e)$ in $\mathds{k}^{\times}$ divides the order of $e$ in $E$.
\end{Lemma}
\begin{proof}
By restricting $(\omega,\beta)$ to $E$, we may assume that $A=E$. In one direction, let us assume that the order of $q(a)$ divides the order of $a$. Now, let us write $$A:=\bigoplus_{k=1}^K\mathbb{Z}/n_k,$$ with $n_k\geq 2$ positive integers, and let $e_k$ denote the generator $1$ in the $k$-th summand. Quinn's formula given in subsection 2.5.2 of \cite{Qui} provides us with an explicit abelian 3-cocycle $(\omega',\beta')$ associated to the quadratic form $q$. In particular, we have $$\omega'(a,b,c):=\prod^K_{\substack{k=1\\b_k+c_k>n_k}}q(e_k)^{a_kn_k},$$ for every $a,b,c\in A$. But, by hypothesis, we have $\omega'= 1$, and by \cite{ML}, $(\omega,\beta)$ and $(\omega',\beta')$ differ by a coboundary. Thence, $\omega$ is trivializable.

Conversely, let us assume that there exists $a\in A$ such that the order of $q(a)$ is strictly greater than $\mathrm{ord}(a)$. As $q$ is a quadratic form, the order of $q(a)$ divides $\mathrm{gcd}(\mathrm{ord}(a)^2,2\cdot \mathrm{ord}(a))$, thus $\mathrm{ord}(a)=2n$ is even, and $q(a)^{\mathrm{ord}(a)}=-1$. Now, using Quinn's formula again, the restriction of $\omega$ to $\langle n\cdot a\rangle\cong\mathbb{Z}/2$ is cohomologous to the 3-cocycle given by $$\omega''(k,l,m):=\begin{cases} -1,& \mathrm{if\ }k,l,m\mathrm{\ are\ odd}\\ 1,& \mathrm{otherwise}\end{cases},$$ for every $k,l,m\in\mathbb{Z}/2$. But, $\omega''$ is non-trivial in cohomology, so $\omega$ is non-trivial in cohomology. This finishes the proof of the result.
\end{proof}

\begin{Corollary}
If $E$ is a subgroup of $A$ of odd order, then the restriction of $\omega$ to $E$ is trivializable.
\end{Corollary}
\begin{proof}
As $q$ is a quadratic form, we know that for every $a\in A$ the order of $q(a)$ divides $\mathrm{gcd}(\mathrm{ord}(a)^2,2\cdot \mathrm{ord}(a))$, where $\mathrm{ord}(a)$ is the order of $a$. Now, if $E$ has odd order, then every object $e\in E$ has odd order, so $\mathrm{gcd}(\mathrm{ord}(e)^2,2\mathrm{ord}(e))=\mathrm{ord}(e)$. This implies that $q(e)^{\mathrm{ord}(e)}=1$ for every $e\in E$. The result follows from lemma \ref{lem:subgrouptrivializable}.
\end{proof}

\begin{Remark}
It is tempting to believe that there is a unique largest subgroup of $A$, over which the restriction of $\omega$ is trivializable. However, this is not the case in general, as can be seen from example \ref{ex:Z/4+Z/2} below.
\end{Remark}

\subsection{A Decomposition Principle}

We now wish to prove that the relative Deligne tensor product over a pointed braided fusion category can be computed ``one prime at a time''. This justifies why all of the examples we will consider later are finite abelian $p$-groups for some prime $p$. For completeness, we begin by showing that every pointed braided fusion category can be split into a Deligne tensor product of pointed braided fusion categories whose groups of equivalence classes of simple objects have coprime order.

\begin{Lemma}\label{lem:quadraticformdecomposition}
Let $A$ and $B$ be two finite abelian groups whose orders are coprime. Any pointed braided fusion category whose group of equivalence classes of invertible object is $A\times B$ is equivalent to $\mathbf{Vect}_A^{q_1}\boxtimes\mathbf{Vect}_B^{q_2}$ for some quadratic forms $q_1:A\rightarrow \mathds{k}^{\times}$ and $q_2:B\rightarrow \mathds{k}^{\times}$.
\end{Lemma}
\begin{proof}
By \cite{ML}, quadratic forms $C\rightarrow \mathds{k}^{\times}$ on a finite abelian group $C$ are in bijective correspondence with $H^3_{ab}(C;\mathds{k}^{\times})$. Further, this last group is by definition isomorphic to the singular cohomology group $H^4(K(C,2);\mathds{k}^{\times})$ of $K(C,2)$, the second Eilenberg-MacLane space associated to the abelian group $C$. This space (or more precisely homotopy type), is uniquely characterized by the property that $\pi_n(K(C,2))=C$ if $n=2$, and $0$ otherwise. In particular, setting $C=A\times B$, we have that $K(A\times B,2)\simeq K(A,2)\times K(B,2)$. But the orders of $A$ and $B$ are coprime, so that, for all positive integers $m$ and $n$, the abelian groups $H_m(K(A,2);\mathbb{Z})$ and $H_n(K(B,2);\mathbb{Z})$ are finite of coprime order by \cite{S1}. Then, we find that $$H_n(K(A\times B,2);\mathbb{Z})\cong H_n(K(A,2);\mathbb{Z})\times H_n(K(B,2);\mathbb{Z})$$ for all $n$ by applying the K\"unneth formula. Further, using the universal coefficient theorem for $\mathds{k}^{\times}$, we get $$H^4(K(A\times B,2);\mathds{k}^{\times})\cong H^4(K(A,2);\mathds{k}^{\times})\times H^4(K(B,2);\mathds{k}^{\times}).$$ Inspection shows that this isomorphism has to be implement by the canonical maps, which concludes the proof.
\end{proof}

We also need the following result which follows from the first part of the proof of theorem 4.6 of \cite{D3} by recalling that the Frobenius-Perron dimension of a pointed fusion category is equal to the order of its group of equivalence classes of invertible objects.

\begin{Proposition}
Let $\mathcal{C}$ and $\mathcal{D}$ be two pointed fusion categories whose groups of equivalence classes of invertible objects have coprime order. Then, for any indecomposable finite semisimple right $\mathcal{C}\boxtimes\mathcal{D}$-module category $\mathcal{P}$, there exists an indecomposable finite semisimple right $\mathcal{C}$-module category $\mathcal{M}$ and an indecomposable finite semisimple right $\mathcal{D}$-module category $\mathcal{N}$ both unique up to equivalence such that $\mathcal{M}\boxtimes\mathcal{N}\simeq \mathcal{P}$ as right $\mathcal{C}\boxtimes\mathcal{D}$-module categories.
\end{Proposition}

By combining the above result together with proposition 3.8 of \cite{DSPS14}, we obtain the desired corollary.

\begin{Corollary}\label{cor:primedecomposition}
Let $\mathcal{C}$, $\mathcal{D}$ be two pointed braided fusion categories whose groups of equivalence classes of invertible objects have coprime order. Then, the action of $\boxtimes_{\mathcal{C}\boxtimes\mathcal{D}}$ on the indecomposable finite semisimple $\mathcal{C}\boxtimes\mathcal{D}$-module categories is completely determined by the actions of $\boxtimes_{\mathcal{C}}$ 
on the indecomposable finite semisimple $\mathcal{C}$-module categories and of $\boxtimes_{\mathcal{D}}$ on the finite semisimple $\mathcal{D}$-module categories
\end{Corollary}

\section{Main Results}

In order to describe the relative Deligne tensor product of the indecomposable finite semisimple module categories over pointed braided fusion categories, we first need to prove a technical result.

\subsection{The Technical Heart}

Let us fix $A$ a finite abelian group, and $\omega$ a 3-cocycle with coefficient in $\mathds{k}^{\times}$, which we will assume to be normalized. Given $B$ a finite abelian group, $f:B\rightarrow A$ a homomorphism, and $\phi:B\times B\rightarrow \mathds{k}^{\times}$ be a 2-cochain such that $d\phi=f^*\omega^{-1}$, we can consider the algebra $\mathds{k}^{\phi}\lbrack B\rbrack$ in $\mathbf{Vect}^{\omega}_A$. We wish to describe the category of module over this algebra. In order to do so, let us write $K=Ker(f)$.

\begin{Lemma}\label{lem:bilinear}
The function $Bil(\phi):B\times K\rightarrow \mathds{k}^{\times}$ given by $$Bil(\phi)(b,k):=\phi(b,k)/\phi(k,b)$$ is bilinear.
\end{Lemma}
\begin{proof}
Let us begin by observing that, as $\omega$ is normalized, for every $b,c,d$ in $B$, $\omega(f(b),f(c),f(d))=1$ whenever $a,b,$ or $c$ is in $K$. Thus, by repeated application of the 2-cocycle condition for $\phi$, we have $$\frac{\phi(b,k)\phi(c,k)}{\phi(k,b)\phi(k,c)}=\frac{\phi(b+c,k)\phi(k+b,c)\phi(b,k)}{\phi(k,b+c)\phi(b,c+k)\phi(k,c)}=\frac{\phi(b+c,k)}{\phi(k,b+c)},$$ for every $b,c$ in $B$, and $k\in K$. This proves linearity in the first variable. Linearity in the second variable follows using a similar argument.
\end{proof}

Let us set $K^{\bot}:=\{b\in B| Bil(\phi)(b,k)=1\mathrm{\ for\ all\ }k\in K\}$. By means of lemma \ref{lem:bilinear}, we find that $K^{\bot}$ is in fact a subgroup of $B$. Let us also define the subgroup $(K\cap K^{\bot})^{\bot}:=\{b\in B| Bil(\phi)(b,k)=1\mathrm{\ for\ all\ }k\in K\cap K^{\bot}\}$. The following lemma will be used in the proof of proposition \ref{prop:decompositionalgebraicallyclosed}.

\begin{Lemma}\label{lem:orthogonalcomplement}
We have $(K\cap K^{\bot})^{\bot} = K + K^{\bot}$.
\end{Lemma}
\begin{proof}
Let $p=char(\mathds{k})$, and set $B^{\bot}:=\{k\in K|Bil(\phi)(b,k)=1\mathrm{\ for\ all\ }b\in K\}$. Note that the $p$-torsion subgroup $K_p$ of $K$ is contained in $B^{\bot}$, and that the $p$-torsion subgroup of $B$ is contained in $K^{\bot}$. Now, $Bil(\phi)$ induces a linear map $B\rightarrow Hom(K,\mathds{k}^{\bot})$, whose kernel is $K^{\bot}$, so the induced map $B/K^{\bot}\hookrightarrow Hom(K,\mathds{k}^{\bot})$ is injective. Thus, we get a surjection $$Hom(Hom(K,\mathds{k}^{\times}),\mathds{k}^{\times})\twoheadrightarrow Hom(B/K^{\bot},\mathds{k}^{\times}).$$ But the left hand-side is canonically identified with $K/K_p$, and the surjection above may be identified with the map $$K/K_p\rightarrow Hom(B/K^{\bot},\mathds{k}^{\times})$$ induced by $Bil(\phi)$. As the kernel of this last map is $B^{\bot}/K_p$, we get an isomorphism $K/B^{\bot}\cong Hom(B/K^{\bot},\mathds{k}^{\times})$. Further, as $B/K^{\bot}$ has no $p$-torsion, there is a non-canonical isomorphism $Hom(B/K^{\bot},\mathds{k}^{\times})\cong B/K^{\bot}$, so we find $|K/B^{\bot}|=|B/K^{\bot}|$.

Using an analogous reasoning, we find that $|B/(K\cap K^{\bot})^{\bot}| = |(K\cap K^{\bot})/B^{\bot}|$. Furthermore, note that $K+K^{\bot}$ is certainly contained in $(K\cap K^{\bot})^{\bot}$. Combining this fact with the equalities derived above, we find $$|K+K^{\bot}|\leq |(K\cap K^{\bot})^{\bot}| = \frac{|B|\cdot |B^{\bot}|}{|K\cap K^{\bot}|}=\frac{|B|\cdot |B^{\bot}|\cdot |K+K^{\bot}|}{|K|\cdot |K^{\bot}|}=|K+ K^{\bot}|.$$ This finishes the proof of the lemma.
\end{proof}

Let us denote by $H$ the image of $K^{\bot}$ under $f$. Fixing a set-theoretic section $s:H\hookrightarrow K^{\bot}$ of the projection $K^{\bot}\twoheadrightarrow H$, we also define a 2-cochain $\rho:H\times H\rightarrow \mathds{k}^{\times}$ by \begin{equation}\label{eqn:rhoprop}\rho(g,h):=\omega^{-1}(-h-g,g,h)\omega(-h,-g,g)\phi^{-1}(-s(h),-s(g))\end{equation} for every $g,h\in H$.

\begin{Proposition}\label{prop:decompositionalgebraicallyclosed}
The finite semisimple right $\mathbf{Vect}^{\omega}_A$-module category $\mathcal{N}$ of left $\mathds{k}^{\phi}\lbrack B\rbrack$-modules in $\mathbf{Vect}^{\omega}_A$ is equivalent to $\oplus_{i=1}^m\mathcal{M}(H,\rho)$ with $$m:=\frac{|K|\cdot|K^{\bot}|}{|B|}.$$
\end{Proposition}
\begin{proof}
We begin by examining the simple objects of $\mathcal{N}$. Given a pair $(a,M)$ consisting of an element $a\in A$ and a simple module $M$ over the algebra $\mathds{k}^{\phi}\lbrack K\rbrack$ in $\mathbf{Vect}$, we can consider $$N(a,M):=\mathds{k}^{\phi}\lbrack B\rbrack\otimes_{\mathds{k}^{\phi}\lbrack K\rbrack} (M\otimes \mathds{k}_a),$$ which is a simple object of $\mathcal{N}$. Moreover, every simple object of $\mathcal{N}$ is of this form. Namely, given any simple object $N$ of $\mathcal{N}$, there exists $a\in A$ such that $N_a$ the $a$-graded part of $N$ is non-zero. As $N$ is simple, $N_a$ must be a simple $\mathds{k}^{\phi}\lbrack K\rbrack$-module, and we must have $N\cong \mathds{k}^{\phi}\lbrack B\rbrack\otimes_{\mathds{k}^{\phi}\lbrack K\rbrack} N_a.$ Now, there are some redundancies in this description. More precisely, given $b\in B$, then the associator provides an isomorphism $N(a,M)\otimes \mathds{k}_{f(b)} \cong N(a+f(b),M)$. Further, a quick computation shows that the $a$-graded part of $N(a+f(b),M)$ is \begin{equation}\label{eqn:decompositioniso}\begin{tabular}{rccc} $\Sigma^b:$&$M^{\phi(b)}\otimes \mathds{k}_a$&$\xrightarrow{\simeq}$&$\mathds{k}^{\phi}\lbrack -b + K \rbrack \otimes_{\mathds{k}^{\phi}\lbrack K\rbrack} (M\otimes \mathds{k}_{a+f(b)}),$\\&$m\otimes 1_a$&$\mapsto$&$(-b)\otimes m\otimes 1_{a+f(b)}$\end{tabular}\end{equation} where $M^{\phi(b)}$ denotes $M$ with the action twisted by multiplying with $Bil(\phi)(b,-)$. This means that we have $N(a,M)\otimes \mathds{k}_{f(b)} \cong N(a,M^{\phi(b)})$. We emphasize that the isomorphism $\Sigma_h$ does depend on $b$, and not merely on $f(b)$.

Observe that given two simple $\mathds{k}^{\phi}\lbrack K\rbrack$-modules $M$ and $M'$, $N(a,M)$ and $N(a,M')$ are isomorphic as $\mathds{k}^{\phi}\lbrack B\rbrack$-modules if and only if $M\cong M'$ as $\mathds{k}^{\phi}\lbrack K\rbrack$-modules. Combining this observation with the results of the previous paragraph, we find that simple objects of $\mathcal{N}$ up to equivalence are parametrised by pairs consisting of a class in $A/f(B)$ and a simple $\mathds{k}^{\phi}\lbrack K\rbrack$-module. As $\mathds{k}^{\phi}\lbrack K\rbrack$ is a finite semisimple algebra with center $\mathds{k}\lbrack K\cap K^{\bot}\rbrack$, $\mathcal{N}$ has $|K\cap K^{\bot}|\cdot |A/f(B)|$ equivalence classes of simple objects. But, for any 2-cochains $\rho_i:H\times H\rightarrow \mathds{k}^{\times}$ such that $d\rho_i=\omega|_H^{-1}$, the finite semisimple category $\mathcal{M}(H,\rho_i)$ has $|A/H|$ many equivalence classes of simple objects. So we find that $\oplus_{i=1}^m \mathcal{M}(H,\rho_i)$ has $$m\cdot |A/H| = \frac{|K|\cdot|K^{\bot}|\cdot |A|}{|B|\cdot |H|}=\frac{|K^{\bot}|}{|H|}\cdot \frac{|K|\cdot |A|}{|B|}=|K\cap K^{\bot}|\cdot |A/f(B)|$$ equivalence classes of simple objects, i.e. $\oplus_{i=1}^m \mathcal{M}(H,\rho_i)$ and $\mathcal{N}$ have the same number of equivalence classes of simple objects.

Let us fix a pair $(a,M)$, with $a\in A$ and $M$ a simple $\mathds{k}^{\phi}\lbrack K\rbrack$-module. The stabilizer $S\subseteq A$ of $N(a,M)$ under the right action by the simple objects of $\mathbf{Vect}^{\omega}_A$ is evidently contained in $f(B)$ as $S$ must preserve the class of $a$ in $A/f(B)$. So let $b\in B$, and consider the right action of $\mathds{k}_{f(b)}$ on $N$. Thanks to our computations above, we know that $N(a,M)\otimes \mathds{k}_{f(b)} \cong N(a,M^{\phi(b)})$. This implies that $f(b)\in S$ if and only if $M=M^{\phi(b)}$ as $\mathds{k}^{\phi}\lbrack K\rbrack$-modules. To check this, it is enough to compare how the center $\mathds{k}\lbrack K\cap K^{\bot}\rbrack$ acts. Thus, $f(b)\in S$ if and only if $Alt(\phi)(b,k)=1$ for every $k\in K\cap K^{\bot}$. But, by lemma \ref{lem:orthogonalcomplement}, we have $(K\cap K^{\bot})^{\bot}=K+K^{\bot}$, so $S=H$. This proves that as right $\mathbf{Vect}_A^{\omega}$-module categories, $\mathcal{N}\simeq \oplus_{i=1}^m \mathcal{M}(H,\rho_i)$ for some 2-cochains $\rho_i:H\times H\rightarrow \mathds{k}^{\times}$ such that $d\rho_i=\omega|_H^{-1}$.

It remains to determine the 2-cochains $\rho_i$. Observe that any of the indecomposable summands of $\mathcal{N}$ contains a simple object of the form $N(0,M)$ for some simple $\mathds{k}^{\phi}\lbrack K\rbrack$-module $M$. Thence, following the proof of lemma \ref{lem:algclosedmodules}, we need to determine the coherence isomorphisms for $N(0,M)$ under the right action by $\mathds{k}_h$ for $h\in H$. But, $\mathcal{N}$ is not skeletal, so we need to pick an isomorphism $N(0,M)\cong N(0,M)\otimes \mathds{k}_h$ for every $h\in H$. In order to do so, we use $s:H\hookrightarrow K\cap K^{\bot}$ the set-theoretic section of $f:K\cap K^{\bot}\twoheadrightarrow H$ fixed above. We set $T_h:N(0,M)\cong N(0,M)\otimes \mathds{k}_h$ to be the isomorphism whose $0$-graded component is $\Sigma_{s(h)}$ as in (\ref{eqn:decompositioniso}). In particular, we have $$T_h(b\otimes m \otimes 1_0)=\omega^{-1}(b,-h,h) \phi(b,-s(h))\cdot \big( (b-s(h))\otimes m \otimes 1_0 \big) \otimes 1_h$$ for every $b\in B$ and $m\in M$. Thus, given any $g,h\in H$, we find that the $0$-graded part of the composite

$$\begin{tikzcd}
{N(0,M)\otimes (\mathds{k}_{g} \otimes \mathds{k}_{h})} \arrow[rr, "\omega^{-1}"] &  & {(N(0,M)\otimes \mathds{k}_{g}) \otimes \mathds{k}_{h}} \arrow[d, "T_g^{-1}\otimes Id"] \\
                                                                                  &  & {N(0,M)\otimes\mathds{k}_h} \arrow[d, "T_h^{-1}"]                                       \\
{N(0,M)} \arrow[uu, "T_{gh}"]                                                     &  & {N(0,M)}                                                                               
\end{tikzcd}$$

\noindent is the automorphism of the simple $\mathds{k}^{\phi}\lbrack K\rbrack$-module $M$ given by multiplication by $\rho(g,h)$, as defined in equation (\ref{eqn:rhoprop}). In particular, $\rho$ must necessarily satisfy $d\rho = \omega|_H^{-1}$ because $\mathcal{N}$ is a right $\mathbf{Vect}^{\omega}_A$-module category. As $M$ was arbitrary, $\rho_i=\rho$ for every $i$, so that $\mathcal{N}\simeq \oplus_{i=1}^m\mathcal{M}(H,\rho)$ as right $\mathbf{Vect}^{\omega}_A$-module categories.
\end{proof}

\begin{Remark}\label{rem:omegatrivial1}
If the 3-cocycle $\omega$ is trivial on the image $Im(f)$ of $f:B\rightarrow A$, then, $\phi$ is a 2-cocycle, which means that the function $Alt(\phi):B\times B\rightarrow \mathds{k}^{\times}$ defined by $Alt(\phi)(b,c) = \phi(b,c)/\phi(c,b)$ for every $b,c$ in $B$ is bilinear. In addition, $Alt(\phi)$ is skew-symmetric, i.e. it satisfies $Alt(\phi)(b,b)=1$ for all $b$ in $B$. Now, the restriction of $Alt(\phi)$ to $K^{\bot}$ descends to a skew-symmetric bilinear form $\sigma$ on $H$. But, as $\omega|_{Im(F)}= 1$, $\rho$ is a 2-cocycle as well, so that $Alt(\rho)$ is a skew-symmetric bilinear form on $H$. Furthermore, it is easy to check that $Alt(\rho) = \sigma$. Finally, by proposition 3.6 of \cite{Tam}, $Alt$ defines an isomorphism between $H^2(H;\mathds{k}^{\times})$ and the group of skew-symmetric bilinear forms on $H$, so $\sigma$ determines $\rho$ uniquely up to coboundary. In particular, if $char(\mathds{k})=0$ and $\omega=1$, we recover proposition 2.3 of \cite{ENO2}.
\end{Remark}

\subsection{The Main Theorem and its Corollaries}

Let us now fix a finite abelian group $A$ and $(\omega, \beta)$ an abelian 3-cocycle for $A$ with coefficient in $\mathds{k}^{\times}$. We now compute the relative Deligne tensor product of all indecomposable finite semisimple right $\mathbf{Vect}^{(\omega,\beta)}_A$-module categories. By restricting to the separable ones characterised in lemma \ref{lem:charsep}, this also computes the fusion rule of $\mathbf{Mod}(\mathbf{Vect}^{(\omega,\beta)}_A)$. 

Let $E, F$ be two subgroups of $A$, and $\phi:E\times E\rightarrow\mathds{k}^{\times}$, $\psi:F\times F\rightarrow \mathds{k}^{\times}$ be two cochains such that $d\phi=\omega|_E^{-1}$ and $d\psi=\omega|_F^{-1}$. We define a 2-cochain $\chi$ on $E\oplus F$ with coefficient in $\mathds{k}^{\times}$ by \begin{align*}\chi((e_1,f_1),(e_2,f_2)):=& \omega(e_1,f_1,e_2+f_2)\omega^{-1}(f_1,e_2,f_2)\beta(f_1,e_2)\omega(e_2,f_1,f_2)\\ & \omega^{-1}(e_1,e_2,f_1+f_2)\phi(e_1,e_2)\phi(f_1,f_2)\end{align*} for every $e_1,e_2\in E$ and $f_1,f_2\in F$. Writing $\pi:E\oplus F\rightarrow A$ for the canonical map given by $\pi(e,f)=e+f$, it is easy to check that $d\chi = \pi^*\omega^{-1}$. Furthermore, the kernel of $\pi$ is $E\cap F$ viewed as a subgroup of $E\oplus F$ via $e\mapsto (e,-e)$, thus we can consider $(E\cap F)^{\bot}$, the orthogonal complement of $E\cap F$ with respect to the bilinear map $Bil(\chi):(E\oplus F)\times (E\cap F)\rightarrow \mathds{k}^{\times}$, as defined in lemma \ref{lem:bilinear}. We write $H$ for the image of $(E\cap F)^{\bot}$ under $\pi$. Finally, fixing a set-theoretic section $s:H\hookrightarrow (E\cap F)^{\bot}$ of the projection $(E\cap F)^{\bot}\twoheadrightarrow H$, we define a 2-cochain $\rho:H\times H\rightarrow \mathds{k}^{\times}$ by $$\rho(g,h):=\omega^{-1}(-h-g,g,h)\omega(-h,-g,g)\chi^{-1}(-s(h),-s(g))$$ for every $g,h\in H$.

\begin{Theorem}\label{thm:fusionalgclosed}
Write $\mathcal{C}=\mathbf{Vect}^{(\omega,\beta)}_A$ for the pointed braided fusion category of $A$-graded vector spaces with associator $\omega$ and braiding $\beta$. Let $E, F$ be two subgroups of $A$, and let $\phi:E\times E\rightarrow \mathds{k}^{\times}$, $\psi:F\times F\rightarrow \mathds{k}^{\times}$ be two 2-cochains such that $d\phi=\omega|_E^{-1}$ and $d\psi=\omega|_F^{-1}$. The relative Deligne tensor product of the two right $\mathcal{C}$-module categories $\mathcal{M}(E,\phi)$ and $\mathcal{M}(F,\psi)$ is given by
$$\mathcal{M}(E,\phi)\boxtimes_{\mathcal{C}}\mathcal{M}(F,\psi) \simeq \bigoplus_{i=1}^m \mathcal{M}(H,\rho),$$
where $H$ and $\rho$ are the subgroup and 2-cochain defined in the previous paragraph, and $$m=\frac{|(E\cap F)^{\bot}|\cdot |E\cap F|}{|E|\cdot |F|}.$$
\end{Theorem}
\begin{proof}
By definition, the relative Deligne tensor product $\mathcal{M}(E,\phi)\boxtimes_{\mathcal{C}}\mathcal{M}(F,\psi)$ is given by the category of left modules over the algebra $R(E,\phi)\otimes R(F,\psi)$ in $\mathcal{C}$. Note that, as recalled in equation (\ref{eqn:tensoralgebras}), the multiplication of this algebra is twisted by the associator $\omega$ and the braiding $\beta$ of $\mathcal{C}$. More precisely, it is given by the 2-cochain $\chi$ on $E\oplus F$. Proposition \ref{prop:decompositionalgebraicallyclosed} now yields the desired result.
\end{proof}

If the 3-cocycle $\omega$ is trivial on $E+F$, then $\phi$, $\psi$, and $\rho$ are 2-cocycles. In particular, $Alt(\chi)$ is a skew-symmetric bilinear form on $E\oplus F$, which satisfies $$Alt(\chi)((e_1,f_1),(e_2,f_2))= Alt(\phi)(e_1,e_2)Alt(\psi)(f_1,f_2)\beta(f_1,e_2)/\beta(f_2,e_1),$$ for every $e_1,e_2\in E$ and $f_1,f_2\in F$. Then, as in remark \ref{rem:omegatrivial1}, the restriction of $Alt(\chi)$ to $(E\cap F)^{\bot}$ descends to a skew-symmetric bilinear form $\sigma$ on $H$. It is easy to check that $Alt(\rho)=\sigma$, so that the cohomology class of $\rho$ is uniquely determined by $\sigma$. In practice, this is a really efficient method to determine $\rho$, so we summarize the above discussion in the result below.

\begin{Corollary}
If the 3-cocycle $\omega$ is trivial on $E+F$, then the skew-symmetric bilinear form $Alt(\chi)$ on $E\oplus F$ descends to a skew-symmetric bilinear form $\sigma$ on $H$, such that $Alt(\rho)=\sigma$. Further, this determines the 2-cocycle $\rho$ on $H$ uniquely up to coboundary.
\end{Corollary}

\begin{Remark}
If $char(\mathds{k})=0$ and $\omega=1$, we recover proposition 2.5.1 of \cite{D2}. If we assume in addition that $\beta=1$, we also recover proposition 3.16 of \cite{ENO}.
\end{Remark}

As another corollary, we can generalize proposition 5.5 of \cite{DN1} to algebraically closed fields of positive characteristic. We also expand their result by exhibiting the inverse explicitly.

\begin{Corollary}\label{cor:invertible}
The indecomposable finite semisimple right $\mathbf{Vect}^{(\omega,\beta)}_A$-module category $\mathcal{M}(E,\phi)$ is invertible if and only if the bilinear form $E\times E\rightarrow \mathds{k}^{\times}$ given by $(e,f)\mapsto \beta(e,f)\phi(f,e)/\phi(e,f)$ is non-degenerate. In that case, its inverse is $\mathcal{M}(E,\phi^{op} \beta)$.
\end{Corollary}
\begin{proof}
Write $\mathcal{C}:=\mathbf{Vect}_A^{(\omega,\beta)}$. Let $F$ be a subgroup of $A$ and $\psi:F\times F\rightarrow \mathds{k}^{\times}$ a 2-cochain such that $d\psi = \omega|_F^{-1}$, and $$\mathcal{M}(E,\phi)\boxtimes_{\mathcal{C}}\mathcal{M}(F,\psi)\simeq \mathcal{M}(\langle 0\rangle, triv),$$ where $triv$ denotes the trivial cochain. To avoid ambiguity, let us write $K:=\{(e,-e)\in E\oplus F|e\in E\cap F\}$. The orthogonal complement $K^{\bot}$ of $K$ under the bilinear map $Bil(\chi):(E\oplus F)\times K\rightarrow \mathds{k}^{\times}$ is contained in $K$. But, $Bil(\chi)$ induces a linear map $\widehat{Bil(\chi)}:F\rightarrow Hom(K,\mathds{k}^{\times})$ by sending $f$ to $e\mapsto Bil(\chi)((f,0),(e,-e))$. But theorem \ref{thm:fusionalgclosed} together with our assumption implies that the kernel of $\widehat{Bil(\chi)}$ is contained in $(\langle 0\rangle\oplus F)\cap K^{\bot}$, which is trivial. This forces $|E\cap F|\geq |F|$, and so $E\cap F = F$. Similarly, one proves that $E=E\cap F$, and so we must have $E=F$.

Now, let $\phi':A\times A\rightarrow \mathds{k}^{\times}$ be the extension by zero of the 2-cochain $\phi:E\times E\rightarrow \mathds{k}^{\times}$. Using $\phi'$ we find that the two abelian 3-cocycles $(\omega,\beta)$ and $(\omega',\beta')$ are equivalent, where $$\omega'=\omega\cdot d\phi'\mathrm{\ and\ }\beta'(a,b):=\beta(a,b)\phi'(b,a)/\phi'(a,b)$$ for all $a,b\in A$. In particular, we get a braided monoidal equivalence $T:\mathbf{Vect}_A^{(\omega,\beta)}\rightarrow \mathbf{Vect}_A^{(\omega',\beta')}$, whose underlying functor is the identity, and with coherence isomorphisms $T_{\mathds{k}_a,\mathds{k}_b}:\mathds{k}_a\otimes\mathds{k}_b\cong\mathds{k}_a\otimes\mathds{k}_b$ is given by multiplication by $\phi'(a,b)^{-1}$ for all $a,b\in A$. Note that, by construction, $\omega'|_E=1$. Further, writing $\mathcal{C}':= \mathbf{Vect}_A^{(\omega',\beta')}$, the image of the algebra $R(E,\phi)$ in $\mathcal{C}$ under $T$ is the algebra $T(R(E,\phi))=R(E,triv)$ in $\mathcal{C}'$, and likewise $T(R(E,\psi))=R(E,\psi\phi^{-1})$ as algebras in $\mathcal{C}'$. Let us write $\mathcal{M}'(E,triv)=Mod_{\mathcal{C}'}(R(E,triv))$, and $\mathcal{M}'(E,\psi\phi^{-1})=Mod_{\mathcal{C}'}(R(E,\psi\phi^{-1}))$. Now, by definition, $\mathcal{M}(E,\phi)\boxtimes_{\mathcal{C}}\mathcal{M}(E,\psi)=Mod_{\mathcal{C}}(R(E,\phi)\otimes R(E,\psi))$. Further, using the braided monoidal equivalence $T$, we find that the latter category is equivalent to $$Mod_{\mathcal{C}'}(T(R(E,\phi)\otimes R(E,\psi)))= \mathcal{M}'(E,triv)\boxtimes_{\mathcal{C}'}\mathcal{M}'(E,\psi\phi^{-1}),$$ and that $\mathcal{M}(E,\phi)\boxtimes_{\mathcal{C}}\mathcal{M}(E,\phi)\simeq \mathcal{M}(\langle 0\rangle,triv)$ if and only if $\mathcal{M}'(E,triv)\boxtimes_{\mathcal{C}'}\mathcal{M}'(E,\psi\phi^{-1})\simeq \mathcal{M}'(\langle 0\rangle,triv)$. Let us now assume that the statement of the corollary holds for $\mathcal{C}'$. Under this hypothesis, we find that the bilinear form $Alt(triv)\cdot \beta' = \beta':E\times E\rightarrow \mathds{k}^{\times}$ is non-degenerate if and only if $\mathcal{M}(E,\phi)$ is invertible. Further, if this is the case, we see that $\psi\phi^{-1}$ is cohomologous to $\beta'$, so that $\psi$ is cohomologous to $\phi^{op}\beta$. It therefore only remains to prove the result under the additional assumption that $\omega|_E=1$.

\textbf{Assuming $\omega$ is trivial on $E$:} Assume that $\mathcal{M}(E,\phi)$ is invertible. By the first paragraph of the proof of this corollary, its inverse is of the form $\mathcal{M}(E,\psi)$, where $\psi:E\times E\rightarrow \mathds{k}^{\times}$ is a 2-cocycle. So, let us fix a 2-cocycle $\psi$ such that $$\mathcal{M}(E,\phi)\boxtimes_{\mathcal{C}}\mathcal{M}(E,\psi)\simeq \mathcal{M}(\langle 0\rangle, triv).$$ Recall that we write $K=\{(e,-e)|e\in E\}\subseteq E\oplus E$. By theorem \ref{thm:fusionalgclosed}, we have $|K^{\bot}|\cdot |K| = |E|\cdot |E|$, and, as seen above, we must have $K^{\bot}\subseteq K$, so we find $K^{\bot}=K$. In particular, for every $e,f$ in $E$, we must have $Alt(\chi)((e,-e),(f,-f)) =1$, but, as $\omega|_E=1$, this gives $$Alt(\phi)(e,f)=Alt(\psi)(e,f)^{-1}\cdot \frac{\beta(e,f)}{\beta(f,e)}$$ for every $e,f$ in $E$. Thus, $\psi$ must differ from $\phi^{op}\beta$ by a coboundary.

Let us finish the proof by showing that $$\mathcal{M}(E,\phi)\boxtimes_{\mathcal{C}}\mathcal{M}(E,\phi^{op} \beta)\simeq \mathcal{M}(\langle 0\rangle, triv)=\mathcal{C}$$ if and only if $Alt(\phi)\cdot \beta$ is non-degenerate. We consider the skew-symmetric bilinear form $\tau$ on $E\oplus E$ given by $$\tau((e_1,e_2),(f_1,f_2)) = \frac{Alt(\phi)(e_1,f_1)}{Alt(\phi)(e_2,f_2)}\cdot \frac{\beta(e_2,f_1)}{\beta(f_2,e_1)}\cdot \frac{\beta(e_2,f_2)}{\beta(f_2,e_2)}.$$ The orthogonal complement of $K$ in $E\oplus E$ contains $K$. But, for any $e,f$ in $E$, we have $$\tau((0,e),(f,-f))=\frac{\beta(e,f)}{Alt(\phi)(e,-f)}=Alt(\phi)(e,f)\beta(e,f),$$ thence $Alt(\phi)\cdot \beta$ is non-degenerate if and only if $(0,e)\notin K^{\bot}$ for every non-zero $e\in E$. As $K\subseteq K^{\bot}$, this holds if and only if $K^{\bot} = K$. Theorem \ref{thm:fusionalgclosed} then concludes the proof of the result.
\end{proof}

\begin{Remark}
In general, the indecomposable finite semisimple right $\mathbf{Vect}^{(\omega,\beta)}_A$-module categories $\mathcal{M}(E,\phi)$ and $\mathcal{M}(E,\phi^{op}\cdot \beta)$ are dual to one another. Namely, writing $\mathcal{C}=\mathbf{Vect}^{(\omega,\beta)}_A$, corollary 2.4.14 of \cite{DSPS13} implies that the dual of $\mathcal{M}(E,\phi)$ is given by $RMod_{\mathcal{C}}(\mathds{k}^{\phi}\lbrack E\rbrack)$ the category of right $\mathds{k}^{\phi}\lbrack E\rbrack$-modules in $\mathcal{C}$. Note that $RMod_{\mathcal{C}}(\mathds{k}^{\phi}\lbrack E\rbrack)$ is a right $\mathcal{C}$-module category as $\mathcal{C}$ is braided. Now, an explicit computation gives $$RMod_{\mathcal{C}}(\mathds{k}^{\phi}\lbrack E\rbrack)\simeq \mathcal{M}(E,\phi^{op}\beta),$$ and so the claim follows.
\end{Remark}

\begin{Corollary}\label{cor:cyclicorder2}
Let $\mathcal{M}(E,\phi)$ be an invertible finite semisimple right $\mathbf{Vect}^{(\omega,\beta)}_A$-module category such that $E$ is cyclic, then $\mathcal{M}(E,\phi)$ has order at most $2$ in $Pic(\mathbf{Vect}^{(\omega,\beta)}_A)$.
\end{Corollary}
\begin{proof}
As $H^2(E,\mathds{k}^{\times})$ is trivial, $\phi^{op}\beta$ and $\phi$ are cohomologous. By corollary \ref{cor:invertible}, this means that $\mathcal{M}(E,\phi)$ is equivalent to its inverse.
\end{proof}

\begin{Remark}
If $\mathcal{M}(E,\phi)$ is an invertible finite semisimple right $\mathbf{Vect}^{(\omega,\beta)}_A$-module category such that $E$ is generated by at least two elements, its order in the Picard group can be arbitrarily large as can be seen from example \ref{ex:PicZ/p+Z/p} below.
\end{Remark}

\section{Examples}

\begin{Example}
Let $A=\mathbb{Z}/p$. We begin by examining the relative Deligne tensor product over the fusion category $\mathcal{C}:=\mathbf{Vect}_{\mathbb{Z}/p}$ equipped with the trivial braiding. The results are recorded in the table below.

\begin{center}
\begin{tabular}{ |c| c c | } 
\hline
$\boxtimes_{\mathcal{C}}$ & $\langle 0 \rangle$ & $\langle 1 \rangle$\\ 
\hline
$\langle 0 \rangle$ & $\langle 0 \rangle$ & $\langle 1 \rangle$\\ 

$\langle 1 \rangle$ & $\langle 1 \rangle$ & $p\langle 1 \rangle$\\ 
\hline
\end{tabular}
\end{center}

As we have recalled above, abelian 3-cocycles on $A$ up to equivalence correspond to quadratic forms on $A$. But a quadratic form $q$ on $A$ is completely determined by $q(1)$, which is a $\mathrm{gcd}(2p,p^2)$-th root of unity. So assuming that $char(\mathds{k})\neq p$, we can fix $\zeta$ a primitive $p$-th root of unity, and consider the quadratic form given by $q(1)=\zeta$. This determines a braiding $\beta$ on the fusion category $\mathbf{Vect}_{\mathbb{Z}/p}$ by $\beta(1,1)=\zeta$. The relative Deligne tensor product of the indecomposable finite semisimple module categories over $\mathcal{D}:=\mathbf{Vect}^{\beta}_{\mathbb{Z}/p}$ are given below.

\begin{center}
\begin{tabular}{ |c| c c | } 
 \hline
 $\boxtimes_{\mathcal{D}}$ & $\langle 0 \rangle$ & $\langle 1 \rangle$\\ 
 \hline
 $\langle 0 \rangle$ & $\langle 0 \rangle$ & $\langle 1 \rangle$\\ 
 
 $\langle 1 \rangle$ & $\langle 1 \rangle$ & $\langle 0 \rangle$\\ 
 \hline
\end{tabular}
\end{center}

Finally, if $p=2\neq char(\mathds{k})$, there are two additional quadratic forms on $A$ whose value $q(1)$ is a primitive fourth root of unity. In these last cases, the associator of the corresponding fusion category is cohomologically non-trivial. In particular, there is only one indecomposable module category, the trivial one, which is the unit for the relative Deligne tensor product.
\end{Example}

\begin{Example}\label{ex:algclosedZ/p^2}
Let $A=\mathbb{Z}/p^2$. A quadratic form $q$ on $\mathbb{Z}/p^2$ is completely determined by $q(1)$, which has to be a $\mathrm{gcd}(2p^2,p^4)$-th root of unity. Let $q_0:\mathbb{Z}/p^2\rightarrow \mathds{k}^{\times}$ be the quadratic form given by $q_0(1)=1$. We begin by considering the relative Deligne tensor product over $\mathcal{C}:=\mathbf{Vect}_{\mathbb{Z}/p^2}^{q_0}$. The computations of the relative Deligne tensor products over $\mathcal{C}$ are given in the following table.

\begin{center}
\begin{tabular}{ |c| c c c| } 
 \hline
 $\boxtimes_{\mathcal{C}}$ & $\langle 0 \rangle$ & $\langle p \rangle$ & $\langle 1 \rangle$\\ 
 \hline
 $\langle 0 \rangle$ & $\langle 0 \rangle$ & $\langle p \rangle$ & $\langle 1 \rangle$\\ 
 $\langle p \rangle$ & $\langle p \rangle$ & $p\langle p \rangle$ & $p\langle 1 \rangle$\\ 
 $\langle 1 \rangle$ & $\langle 1 \rangle$ & $p\langle 1 \rangle$ & $p^2\langle 1 \rangle$\\ 
 \hline
\end{tabular}
\end{center}

If $p\neq char(\mathds{k})$, the field $\mathds{k}$ has $p^2$-th root of unity, and so there are non-trivial quadratic forms on $\mathbb{Z}/p^2$. Namely, let $\zeta$ be a primitive $p^2$-th root of unity. We let $q_1,q_2:\mathbb{Z}/p^2\rightarrow \mathds{k}^{\times}$ be the quadratic forms given by $q_1(1)=\zeta^p$, and $q_2(1)=\zeta$. We write $\mathcal{D}=\mathbf{Vect}_{\mathbb{Z}/p^2}^{q_1}$ and $\mathcal{E}=\mathbf{Vect}_{\mathbb{Z}/p^2}^{q_2}$ for the corresponding pointed braided fusion categories. In these cases, the associator of the underlying fusion category is trivial, so $\mathcal{D}$ and $\mathcal{E}$ have three non-equivalent indecomposable finite semisimple module categories. The relative Deligne tensor product over $\mathcal{D}$ and $\mathcal{E}$ are recorded in the tables below.

\begin{center}
\begin{tabular}{c c c c}
\begin{tabular}{ |c| c c c| } 
 \hline
 $\boxtimes_{\mathcal{D}}$ & $\langle 0 \rangle$ & $\langle p \rangle$ & $\langle 1 \rangle$\\ 
 \hline
 $\langle 0 \rangle$ & $\langle 0 \rangle$ & $\langle p \rangle$ & $\langle 1 \rangle$\\ 
 $\langle p \rangle$ & $\langle p \rangle$ & $p\langle p \rangle$ & $p\langle 1 \rangle$\\ 
 $\langle 1 \rangle$ & $\langle 1 \rangle$ & $p\langle 1 \rangle$ & $p\langle p \rangle$\\ 
 \hline
\end{tabular}  & & & \begin{tabular}{ |c| c c c| } 
 \hline
 $\boxtimes_{\mathcal{E}}$ & $\langle 0 \rangle$ & $\langle p \rangle$ & $\langle 1 \rangle$\\ 
 \hline
 $\langle 0 \rangle$ & $\langle 0 \rangle$ & $\langle p \rangle$ & $\langle 1 \rangle$\\ 
 $\langle p \rangle$ & $\langle p \rangle$ & $p\langle p \rangle$ & $\langle p \rangle$\\ 
 $\langle 1 \rangle$ & $\langle 1 \rangle$ & $\langle p \rangle$ & $\langle 0 \rangle$\\ 
 \hline
\end{tabular}
\end{tabular}
\end{center}

\noindent In particular, as can be seen from the tables above, the braiding on $\mathbf{Vect}_{\mathbb{Z}/p^2}$ has a big impact on the relative Deligne tensor product!

It remains to investigate the case $p=2\neq char(\mathds{k})$. In that case, the value $q(1)$ of a quadratic form on $\mathbb{Z}/4$ is an $8$-th root of unity. We have already examined the situation when $q(1)$ is a fourth root of unity, so let $q_3$ be a quadratic form on $\mathbb{Z}/4$ such that $q_3(1)$ is a primitive $8$-th root of unity. We write $\mathcal{F}:=\mathbf{Vect}_{\mathbb{Z}/4}^{q_3}$ for the associated pointed braided fusion category. The associator $\omega$ of $\mathcal{F}$ is the element of order $2$ in $H^3(\mathbb{Z}/4;\mathds{k}^{\times})$, and so is trivial on $\mathbb{Z}/2\subseteq \mathbb{Z}/4$. This means that $\mathcal{F}$ has two non-equivalent indecomposable finite semisimple module categories. The computation of the relative Deligne tensor product over $\mathcal{F}$ is given in the following table.

\begin{center}
\begin{tabular}{ |c| c c| } 
 \hline
 $\boxtimes_{\mathcal{F}}$ & $\langle 0 \rangle$ & $\langle 2 \rangle$\\ 
 \hline
 $\langle 0 \rangle$ & $\langle 0 \rangle$ & $\langle 2 \rangle$\\ 
 $\langle 2 \rangle$ & $\langle 2 \rangle$ & $\langle 0 \rangle$\\ 
 \hline
\end{tabular}
\end{center}

\end{Example}

\begin{Example}\label{ex:algclosedZ2Z2}
Let $A=\mathbb{Z}/2\oplus \mathbb{Z}/2$, and write $a=(1,0)$, and $b=(0,1)$ in $\mathbb{Z}/2\oplus \mathbb{Z}/2$. For now, assume $char(\mathds{k})=2$, and let $\mathcal{B}$ be the braided fusion category $\mathbf{Vect}_{\mathbb{Z}/2\oplus\mathbb{Z}/2}$ equipped with the trivial associator and the trivial braiding. Note that these are the only possible choices of associator and braiding. The fusion category $\mathbf{Vect}_{\mathbb{Z}/2\oplus\mathbb{Z}/2}$ has 5 non-equivalent indecomposable finite semisimple module categories corresponding to the 5 subgroups of $\mathbb{Z}/2\oplus \mathbb{Z}/2$ as all the relevant cohomology groups vanish. The relative Deligne tensor products of these indecomposable finite semisimple right $\mathcal{B}$-module categories are given in the table below.

\begin{center}
\begin{tabular}{ |c| c c c c c| } 
 \hline
 $\boxtimes_{\mathcal{B}}$ & $\langle 0 \rangle$ & $\langle a \rangle$ & $\langle b \rangle$ & $\langle a+b \rangle$ & $\langle a, b \rangle$\\ 
 \hline
 $\langle 0 \rangle$ & $\langle 0 \rangle$ & $\langle a \rangle$ & $\langle b \rangle$ & $\langle a+b \rangle$ & $\langle a,b \rangle$\\ 
 $\langle a \rangle$ & $\langle a \rangle$ & $2\langle a \rangle$ & $\langle a,b \rangle$ & $\langle a,b \rangle$ & $2\langle a,b \rangle$\\ 
 $\langle b \rangle$ & $\langle b \rangle$ & $\langle a,b \rangle$ & $2\langle b \rangle$ & $\langle a,b \rangle$ & $2\langle a,b \rangle$\\ 
 $\langle a+b \rangle$ & $\langle a+b \rangle$ & $\langle a,b \rangle$ & $\langle a,b \rangle$ & $2\langle a+b \rangle$ & $2\langle a,b \rangle$\\ 
 $\langle a,b \rangle$ & $\langle a,b \rangle$ & $2\langle a,b \rangle$ & $2\langle a,b \rangle$ & $2\langle a,b \rangle$ & $4\langle a,b \rangle$\\ 
 \hline
\end{tabular}
\end{center}

Now, let us assume that $char(\mathds{k})\neq 2$. Any quadratic form $q:\mathbb{Z}/2\oplus \mathbb{Z}/2\rightarrow \mathds{k}^{\times}$ is completely determined by $q(a)$, $q(b)$, and $q(a+b)$, all of which are fourth roots of unity. Moreover, we must have $q(a)^2q(b)^2q(a+b)^2=1$. This means that there are 32 quadratic forms $\mathbb{Z}/2\oplus \mathbb{Z}/2\rightarrow \mathds{k}^{\times}$. Let us fix $\zeta$ a primitive fourth root of unity. We will examine six pointed braided fusion categories corresponding to the six quadratic forms described in the table below. In fact, as there are six equivalence classes of quadratic forms on $\mathbb{Z}/2\oplus \mathbb{Z}/2$ under outer automorphisms, the quadratic form below represent all possible non-equivalent braidings $\mathbf{Vect}_{\mathbb{Z}/2\oplus\mathbb{Z}/2}$.

\begin{center}
\begin{tabular}{ |c| c c c| } 
 \hline
  & $a$ & $b$ & $a+b$\\ 
 \hline
 $q_0$ & $1$ & $1$ & $1$\\ 
 $q_1$ & $-1$ & $1$ & $1$\\ 
 $q_2$ & $-1$ & $-1$ & $1$\\ 
 $q_3$ & $-1$ & $-1$ & $-1$\\ 
 $q_4$ & $1$ & $\zeta$ & $\zeta$\\ 
 $q_5$ & $-1$ & $\zeta$ & $\zeta$\\ 
 \hline
\end{tabular}
\end{center}

On one hand, observe that the abelian 3-cocycles produced by Quinn's formula using the quadratic forms $q_0$, $q_1$, and $q_2$ have trivial associators. In this case, the underlying fusion category $\mathbf{Vect}_{\mathbb{Z}/2\oplus \mathbb{Z}/2}$ has 6 non-equivalent indecomposable finite semisimple module categories. This is because $H^2(\mathbb{Z}/2\oplus\mathbb{Z}/2,\mathds{k}^{\times})\cong \mathbb{Z}/2$ as $char(\mathds{k})\neq 2$. In particular, there is a non-trivial element $\nu$ in $H^2(\mathbb{Z}/2\oplus\mathbb{Z}/2,\mathds{k}^{\times})$. By abuse of notation, we also write $\nu$ for the pair $(\mathbb{Z}/2\oplus\mathbb{Z}/2,\nu)$. The relative Deligne tensor product of the indecomposable module categories over $\mathcal{C}:=\mathbf{Vect}^{q_0}_{\mathbb{Z}/2\oplus \mathbb{Z}/2}$, $\mathcal{D}:=\mathbf{Vect}^{q_1}_{\mathbb{Z}/2\oplus \mathbb{Z}/2}$, $\mathcal{E}:=\mathbf{Vect}^{q_2}_{\mathbb{Z}/2\oplus \mathbb{Z}/2}$, and $\mathcal{F}:=\mathbf{Vect}^{q_3}_{\mathbb{Z}/2\oplus \mathbb{Z}/2}$ are given in the tables below.

\begin{center}
\begin{tabular}{ |c| c c c c c c| } 
 \hline
 $\boxtimes_{\mathcal{C}}$ & $\langle 0 \rangle$ & $\langle a \rangle$ & $\langle b \rangle$ & $\langle a+b \rangle$ & $\langle a, b \rangle$ & $\nu$\\ 
 \hline
 $\langle 0 \rangle$ & $\langle 0 \rangle$ & $\langle a \rangle$ & $\langle b \rangle$ & $\langle a+b \rangle$ & $\langle a,b \rangle$ & $\nu$\\ 
 $\langle a \rangle$ & $\langle a \rangle$ & $2\langle a \rangle$ & $\langle a,b \rangle$ & $\langle a,b \rangle$ & $2\langle a,b \rangle$ & $\langle a \rangle$\\ 
 $\langle b \rangle$ & $\langle b \rangle$ & $\langle a,b \rangle$ & $2\langle b \rangle$ & $\langle a,b \rangle$ & $2\langle a,b \rangle$ & $\langle b \rangle$\\ 
 $\langle a+b \rangle$ & $\langle a+b \rangle$ & $\langle a,b \rangle$ & $\langle a,b \rangle$ & $2\langle a+b \rangle$ & $2\langle a,b \rangle$ & $\langle a+b \rangle$\\ 
 $\langle a,b \rangle$ & $\langle a,b \rangle$ & $2\langle a,b \rangle$ & $2\langle a,b \rangle$ & $2\langle a,b \rangle$ & $4\langle a,b \rangle$ & $\langle a,b \rangle$\\ 
 $\nu$ & $\nu$ & $\langle a\rangle$ & $\langle b \rangle$ & $\langle a+b \rangle$ & $\langle a,b \rangle$ & $\langle 0 \rangle$\\ 
 \hline
\end{tabular}
\end{center}

\begin{center}
\begin{tabular}{ |c| c c c c c c| } 
 \hline
 $\boxtimes_{\mathcal{D}}$ & $\langle 0 \rangle$ & $\langle a \rangle$ & $\langle b \rangle$ & $\langle a+b \rangle$ & $\langle a, b \rangle$ & $\nu$\\ 
 \hline
 $\langle 0 \rangle$ & $\langle 0 \rangle$ & $\langle a \rangle$ & $\langle b \rangle$ & $\langle a+b \rangle$ & $\langle a,b \rangle$ & $\nu$ \\ 
 $\langle a \rangle$ & $\langle a \rangle$ & $\langle 0 \rangle$ & $\langle a,b\rangle$ & $\nu$ & $\langle b\rangle$ & $\langle a+b\rangle$ \\ 
 $\langle b \rangle$ & $\langle b \rangle$ & $\nu$ & $2\langle b \rangle$ & $\nu$ &  $\langle b \rangle$ & $2\nu$ \\ 
 $\langle a+b \rangle$ & $\langle a+b \rangle$ & $\langle a,b\rangle$ & $\langle a,b\rangle$ & $2\langle a+b \rangle$ & $2\langle a,b\rangle$ & $\langle a+b \rangle$ \\ 
 $\langle a,b \rangle$ & $\langle a,b \rangle$ & $\langle a+b \rangle$ & $2\langle a,b\rangle$ & $\langle a+b\rangle$ & $\langle a,b\rangle$ & $2\langle a+b\rangle$  \\ 
 $\nu$ & $\nu$ & $\langle b\rangle$ & $\langle b \rangle$ & $2\nu$ & $2\langle b\rangle$ & $\nu$  \\ 
 \hline
\end{tabular}
\end{center}

\begin{center}
\begin{tabular}{ |c| c c c c c c| } 
 \hline
 $\boxtimes_{\mathcal{E}}$ & $\langle 0 \rangle$ & $\langle a \rangle$ & $\langle b \rangle$ & $\langle a+b \rangle$ & $\langle a, b \rangle$ & $\nu$\\ 
 \hline
 $\langle 0 \rangle$ & $\langle 0 \rangle$ & $\langle a \rangle$ & $\langle b \rangle$ & $\langle a+b \rangle$ & $\langle a,b \rangle$ & $\nu$\\ 
 $\langle a \rangle$ & $\langle a \rangle$ & $\langle 0 \rangle$ & $\langle a,b \rangle$ & $\nu$ & $\langle b \rangle$ & $\langle a+b \rangle$\\ 
 $\langle b \rangle$ & $\langle b \rangle$ & $\langle a,b \rangle$ & $\langle 0 \rangle$ & $\nu$ & $\langle a\rangle$ & $\langle a+b \rangle$\\ 
 $\langle a+b \rangle$ & $\langle a+b \rangle$ & $\nu$ & $\nu$ & $2\langle a+b \rangle$ & $\langle a+b \rangle$ & $2\nu$\\ 
 $\langle a,b \rangle$ & $\langle a,b \rangle$ & $\langle b \rangle$ & $\langle a \rangle$ & $\langle a+b \rangle$ & $\langle 0 \rangle$ & $\nu$\\ 
 $\nu$ & $\nu$ & $\langle a+b \rangle$ & $\langle a+b \rangle$ & $2\nu$ & $\nu$ & $2\langle a+b \rangle$\\ 
 \hline
\end{tabular}
\end{center}

\begin{center}
\begin{tabular}{ |c| c c c c c c| } 
 \hline
 $\boxtimes_{\mathcal{F}}$ & $\langle 0 \rangle$ & $\langle a \rangle$ & $\langle b \rangle$ & $\langle a+b \rangle$ & $\langle a, b \rangle$ & $\nu$\\ 
 \hline
 $\langle 0 \rangle$ & $\langle 0 \rangle$ & $\langle a \rangle$ & $\langle b \rangle$ & $\langle a+b \rangle$ & $\langle a,b \rangle$ & $\nu$\\ 
 $\langle a \rangle$ & $\langle a \rangle$ & $\langle 0 \rangle$ & $\nu$ & $\langle a,b \rangle$ & $\langle a+b \rangle$ & $\langle b \rangle$\\ 
 $\langle b \rangle$ & $\langle b \rangle$ & $\langle a,b \rangle$ & $\langle 0 \rangle$ & $\nu$ & $\langle a\rangle$ & $\langle a+b \rangle$\\ 
 $\langle a+b \rangle$ & $\langle a+b \rangle$ & $\nu$ & $\langle a,b \rangle$ & $\langle 0 \rangle$ & $\langle b \rangle$ & $\langle a \rangle$\\ 
 $\langle a,b \rangle$ & $\langle a,b \rangle$ & $\langle b \rangle$ & $\langle a+b \rangle$ & $\langle a \rangle$ & $\nu$ & $\langle 0 \rangle$\\ 
 $\nu$ & $\nu$ & $\langle a+b \rangle$ & $\langle a \rangle$ & $\langle b \rangle$ & $\langle 0 \rangle$ & $\langle a,b \rangle$\\ 
 \hline
\end{tabular}
\end{center}

\noindent Once again, by examining the tables above, we can see how big of an impact the braiding has on the relative Deligne tensor product. In particular, note that the first two tables are symmetric along the diagonal, whereas the third is not. This is because $\mathcal{C}$ and $\mathcal{E}$ are symmetric, whereas $\mathcal{D}$ and $\mathcal{F}$ are not. For later use, we wish to record that $Pic(\mathcal{D})$, the group of invertible right $\mathcal{D}$-module categories, is isomorphic to $\mathbb{Z}/2\oplus \mathbb{Z}/2$, and that $Pic(\mathcal{E})\cong S_3$ the symmetric group on three elements.

On the other hand, observe that the abelian 3-cocycles corresponding to the quadratic forms $q_4$, and $q_5$ have the same associator $\omega$, which is non-trivial. The associator $\omega$ is trivial on $\langle a\rangle$, and is non-trivial on all the other non-zero subgroups of $\mathbb{Z}/2\oplus \mathbb{Z}/2$. Thus, we find that the fusion category $\mathbf{Vect}^{\omega}_{\mathbb{Z}/2\oplus \mathbb{Z}/2}$ has 2 non-equivalent indecomposable finite semisimple module categories. Writing $\mathcal{G}:=\mathbf{Vect}^{q_4}(\mathbb{Z}/2\oplus \mathbb{Z}/2)$ and $\mathcal{H}:=\mathbf{Vect}^{q_5}(\mathbb{Z}/2\oplus \mathbb{Z}/2)$, the corresponding relative Deligne tensor products are computed below.

\begin{center}
\begin{tabular}{c c c c c}
\begin{tabular}{ |c| c c| } 
\hline
$\boxtimes_{\mathcal{G}}$ & $\langle 0 \rangle$ & $\langle a \rangle$\\ 
\hline
$\langle 0 \rangle$ & $\langle 0 \rangle$ & $\langle a \rangle$\\ 
$\langle a \rangle$ & $\langle a \rangle$ & $2\langle a \rangle$\\ 
\hline
\end{tabular}  & & & & \begin{tabular}{ |c| c c| } 
\hline
$\boxtimes_{\mathcal{H}}$ & $\langle 0 \rangle$ & $\langle a \rangle$\\ 
\hline
$\langle 0 \rangle$ & $\langle 0 \rangle$ & $\langle a \rangle$\\ 
$\langle a \rangle$ & $\langle a \rangle$ & $\langle 0 \rangle$\\ 
\hline
\end{tabular}
\end{tabular}
\end{center}
\end{Example}

\begin{Example}\label{ex:Z/4+Z/2}
Let $A=\mathbb{Z}/4\oplus\mathbb{Z}/2$, and write $a=(1,0)$, and $b=(0,1)$ in $\mathbb{Z}/4\oplus\mathbb{Z}/2$. Further, assume $char(\mathds{k})\neq 2$, and let $\zeta$ be a primitive fourth root of unity. We consider the quadratic form $q:A\rightarrow \mathds{k}^{\times}$ given by $q(0)=q(2a)=q(3a+b)=1$ and $q(a)=q(b)=q(3a)=q(2a+b)=\zeta$, and write $\mathcal{C}:=\mathbf{Vect}^q_{\mathbb{Z}/4\oplus\mathbb{Z}/2}$ for the associated pointed braided fusion category. Thanks to lemma \ref{lem:subgrouptrivializable}, we find that $\mathcal{C}$ has exactly four equivalence classes of indecomposable finite semisimple module categories. The relative Deligne tensor product of these module categories over $\mathcal{C}$ is given in the following table.

$$\begin{tabular}{ |c| c c c c| } 
 \hline
 $\boxtimes_{\mathcal{C}}$ & $\langle 0 \rangle$ & $\langle 2a \rangle$ & $\langle a \rangle$ & $\langle a+b \rangle$\\ 
 \hline
 $\langle 0 \rangle$ & $\langle 0 \rangle$ & $\langle 2a \rangle$ & $\langle a \rangle$ & $\langle a+b \rangle$\\ 
 $\langle 2a \rangle$ & $\langle 2a \rangle$ & $2\langle 2a \rangle$ & $\langle 2a \rangle$ & $2\langle a+b \rangle$\\ 
 $\langle a \rangle$ & $\langle a \rangle$ & $\langle 2a \rangle$ & $\langle 0 \rangle$ & $\langle a+b \rangle$\\ 
 $\langle a+b \rangle$ & $\langle a+b \rangle$ & $2\langle a+b \rangle$ & $\langle a+b \rangle$ & $4\langle a+b \rangle$\\ 
 \hline
\end{tabular}$$
\end{Example}

\begin{Example}\label{ex:PicZ/p+Z/p}
Let $A=\mathbb{Z}/p\oplus \mathbb{Z}/p$, and write $a=(1,0)$, $b=(0,1)$ in $A$. Assume that $char(\mathds{k})\neq p$, and let $\zeta$ be a primitive $p$-th root of unity. Given $k\in\mathbb{Z}/p$, we write $\mathcal{C}_k$ for the pointed braided fusion category of $A$-graded vector spaces with braiding given by $\beta(a,a)=\zeta^k$, $\beta(b,b) = \zeta$, and $\beta(a,b)=\beta(b,a)=1$. One checks that all the $2(p+1)$ indecomposable finite semisimple right $\mathcal{C}_k$-module category are invertible if and only if $-k$ is not a quadratic residue. If $-k$ is a quadratic residue, then there are $2p$ invertible $\mathcal{C}_k$-module categories.

Let us now fix $k$ such that $-k$ is not a quadratic residue, we will prove that $Pic(\mathcal{C}_k)\cong D_{2(p+1)}$, the dihedral group of order $2(p+1)$. By corollary \ref{cor:cyclicorder2}, if $E\subseteq A$ has cardinality $p$, then $\mathcal{M}(E,triv)$ is an element of order $2$ of $Pic(\mathcal{C}_k)$. Let us write $\phi$ for an element in $H^2(A;\mathds{k}^{\times})$ such that $Alt(\phi)(a,b)=\zeta$. Let us consider the subgroup of $Pic(\mathcal{C}_k)$ generated by the elements $\mathcal{M}(A,\phi^m)$ for $m\in\mathbb{Z}/p$. We claim that it is cyclic of order $p+1$. To see this, note that an explicit computation using theorem \ref{thm:fusionalgclosed} shows that $$\mathcal{M}(A,\phi^m)\boxtimes_{\mathcal{C}_k}\mathcal{M}(A,\phi^n)\simeq \mathcal{M}(A,\phi^q),\ \mathrm{with}\ q=\frac{mn-k}{m+n}$$ for every $m,n\in\mathbb{Z}/p$. Thus, viewing $m$ as a variable, we find that the order of $\mathcal{M}(A,\phi^n)$ in $Pic(\mathcal{C})$ coincides with the order of the matrix $$M_n:=\begin{pmatrix} n & -k\\ 1&n\end{pmatrix}$$ as an element of $PGL(2,p)$. But our assumption on $k$ guarantees that there exists an $n$ such that $M_n$ has order $p+1$, so the claim follows.

Finally, an explicit computation gives $$\mathcal{M}(\langle a\rangle, triv)\boxtimes \mathcal{M}(A, \phi^m)\simeq \mathcal{M}(\langle (-b,a)\rangle,triv)\simeq\mathcal{M}(A, \phi^{-m})\boxtimes\mathcal{M}(\langle a\rangle, triv)$$ for every $m\in\mathbb{Z}/p$. This establishes the desired isomorphism $Pic(\mathcal{C}_k)\cong D_{2(p+1)}$ if $-k$ is not a quadratic residue. If $-k$ is a quadratic residue, a similar analysis proves that $Pic(\mathcal{C}_k)\cong D_{2(p-1)}$.
\end{Example}

\begin{Remark}
For some specific groups and braidings considered above, surprisingly similar formulas have been independently derived in \cite{RSS} through high energy physics considerations. More precisely, the table given in example \ref{ex:algclosedZ/p^2} for $\mathcal{E}$ agrees with equation 6.15 \cite{RSS} upon fixing $p=2$. More generally, one can recover the mathematical equivalent of formula 5.13 of \cite{RSS} using theorem \ref{thm:fusionalgclosed}. Further, the table given in example \ref{ex:algclosedZ2Z2} for $\mathcal{D}$ essentially agrees with the computations of section 6.3 of \cite{RSS}. Finally, the fact derived in example \ref{ex:PicZ/p+Z/p} that $Pic(\mathcal{C}_k)\cong D_{2(p-1)}$ when $-k$ is not a quadratic residue appears in section 6.4 of \cite{RSS}. At the moment, we do not have a completely rigorous mathematical justification for the agreement between our computations.
\end{Remark}

\begin{Example}\label{ex:PicZ72Z/2Z/2}
Let $A=\mathbb{Z}/2^{\oplus 3}$, and write $e_1=(1,0,0)$, $e_2=(0,1,0)$, and $e_3=(0,0,1)$ in $\mathbb{Z}/2^{\oplus 3}$. Assume that $char(\mathds{k})\neq 2$, and let $\beta:\mathbb{Z}/2^{\oplus 3}\times \mathbb{Z}/2^{\oplus 3}\rightarrow \mathds{k}^{\times}$ be the bilinear form given by $\beta(e_i,e_j)=-1$ if $j=i,i+1$ and $\beta(e_i,e_j)=1$ otherwise. Setting $\mathcal{C}:=\mathbf{Vect}_{A}^{\beta}$, we will now prove that $Pic(\mathcal{C})$ is isomorphic to $S_4$, the symmetric group on four elements.

We begin by determining the number of simple objects of $\mathbf{Mod}(\mathcal{C})$, which coincides with the number of non-equivalent indecomposable finite semisimple right $\mathcal{C}$-module categories by lemma \ref{lem:charsep}. Using corollary \ref{cor:abelianmodulecategoriesequivalences}, we find that $\mathbf{Mod}(\mathcal{C})$ has $1$ simple object corresponding to the subgroup of $A$ of order $1$, and $7$ simple objects corresponding to subgroups of order $2$, $7\cdot 2$ simple objects corresponding to subgroups of order $4$ together with a class in $H^2(\mathbb{Z}/2^{\oplus 2},\mathds{k}^{\times})\cong\mathbb{Z}/2$, and $8$ simple objects corresponding to $A$ together with a class in $H^2(\mathbb{Z}/2^{\oplus 3},\mathds{k}^{\times})\cong\mathbb{Z}/2^{\oplus 3}$. This means that $\mathbf{Mod}(\mathcal{C})$ has $30$ simple objects.

Now, let us write $E_{12}$ for the subgroup of $A$ spanned by $e_1$ and $e_2$. The restriction of $\beta$ to $E_{12}$ corresponds to the quadratic form $q_2$ of example \ref{ex:algclosedZ2Z2}. Thus, the inclusion $E_{12}\hookrightarrow A$ induces a braided monoidal functor $F_{12}:\mathcal{E}\rightarrow \mathcal{C}$. By proposition 3.3.4 of \cite{D5}, this induces a monoidal functor $\mathbf{Mod}(F_{12}):\mathbf{Mod}(\mathcal{E})\rightarrow \mathbf{Mod}(\mathcal{C})$. Further, as $F_{12}$ is faithful, $\mathbf{Mod}(F_{12})$ is injective on simple objects, thence it induces an injective map $Pic(\mathcal{E})\rightarrow Pic(\mathcal{C})$. But we have identified $Pic(\mathcal{E})\cong S_3$ in example \ref{ex:algclosedZ2Z2}. This means that $\mathbf{Mod}(F_{12})(\langle a\rangle)=\langle e_1\rangle$, and $\mathbf{Mod}(F_{12})(\langle b\rangle)=\langle e_2\rangle$ are invertible objects of $\mathbf{Mod}(\mathcal{C})$ and satisfy $$\langle e_1\rangle^2 = \langle 0\rangle,\ \langle e_2\rangle^2 = \langle 0\rangle,\ (\langle e_1\rangle\cdot\langle e_2\rangle)^3= \langle 0\rangle$$ in $Pic(\mathcal{C})$. Similarly, we can consider the subgroup $E_{23}$ of $E$ spanned by $e_2$ and $e_3$. As the restriction of $\beta$ to $E_{23}$ also corresponds to $q_3$, we get an inclusion $\mathbf{Mod}(F_{23}):\mathbf{Mod}(\mathcal{E})\rightarrow \mathbf{Mod}(\mathcal{C})$. We have $\mathbf{Mod}(F_{23})(\langle a\rangle)=\langle e_2\rangle$, and $\mathbf{Mod}(F_{23})(\langle b\rangle)=\langle e_3\rangle$, so that $\langle e_3\rangle$ is invertible in $\mathbf{Mod}(\mathcal{C})$ and satisfy $$\langle e_3\rangle^2 = \langle 0\rangle,\ (\langle e_2\rangle\cdot\langle e_3\rangle)^3= \langle 0\rangle$$ in $Pic(\mathcal{C})$. To be able to conclude, we also need to consider the subgroup $E_{13}$ of $A$ spanned by $e_1$ and $e_3$. The restriction of $\beta$ to $E_{13}$ corresponds to the quadratic from $q_2$ of example \ref{ex:algclosedZ2Z2}. This gives an inclusion $\mathbf{Mod}(F_{13}):\mathbf{Mod}(\mathcal{D})\rightarrow \mathbf{Mod}(\mathcal{C})$. As $\mathbf{Mod}(F_{13})(\langle a\rangle)=\langle e_1\rangle$, $\mathbf{Mod}(F_{13})(\langle b\rangle)=\langle e_3\rangle$, and $Pic(\mathcal{D})\cong\mathbb{Z}/2\oplus\mathbb{Z}/2$, we find that $$\langle e_1\rangle \cdot \langle e_3\rangle = \langle e_3\rangle \cdot \langle e_1\rangle$$ in $Pic(\mathcal{C})$. Bringing all of this together, we find that $\langle e_1\rangle$, $\langle e_2\rangle$, and $\langle e_3\rangle$ generate a subgroup of $Pic(\mathcal{C})$, which satisfy the relations defining $S_4$, and so $S_4$ is a subgroup of $Pic(\mathcal{C})$. But, $Pic(\mathcal{C})$ has order at most $30$ as it is a subset of the set of simple objects of $\mathbf{Mod}(\mathcal{C})$. This forces $Pic(\mathcal{C})\cong S_4$ as desired.
\end{Example}

\begin{Remark}
Let $A=\mathbb{Z}/2^{\oplus n}$, and write $e_i$ for the generator of the $i$-th summand. Assume that $char(\mathds{k})\neq 2$, and let $\beta:\mathbb{Z}/2^{\oplus n}\times \mathbb{Z}/2^{\oplus 3}\rightarrow \mathds{k}^{\times}$ be the bilinear form given by $\beta(e_i,e_j)=-1$ if $j=i,i+1$ and $\beta(e_i,e_j)=1$ otherwise. Setting $\mathcal{C}:=\mathbf{Vect}_{A}^{\beta}$, the techniques of example \ref{ex:PicZ72Z/2Z/2} can be applied to prove that $S_{n+1}$ is a subgroup of $Pic(\mathcal{C})$. However, for $n\geq 2$, this is a strict subgroup of $Pic(\mathcal{C})$, which is identified with $GL(n,\mathbb{F}_2)$ as can be seen using example 5.16 of \cite{DN1}.
\end{Remark}

\bibliography{bibliography.bib}

\end{document}